\documentclass[11pt]{article}

\usepackage{geometry}
\usepackage[parfill]{parskip}
\usepackage{graphicx}

\usepackage{amssymb}

\usepackage{epstopdf}

\usepackage{amssymb}
\usepackage{amsmath}
\usepackage{enumerate}
\usepackage{amsthm}

\newfont{\bbb}{msbm10 scaled\magstep1}
\newfont{\bbbsub}{msbm10}

\newcommand{\ba}{\begin{array}}
\newcommand{\ea}{\end{array}}
\newcommand{\be}{\begin{equation}}
\newcommand{\ee}{\end{equation}}

\newtheorem{theorem}{Theorem}
\newtheorem{lemma}[theorem]{Lemma}
\newtheorem{corollary}[theorem]{Corollary}
\newtheorem{rem}[theorem]{Remark}
\let\oldrem\rem
\renewcommand{\rem}{\oldrem\normalfont}
\newtheorem{definition}[theorem]{Definition}
\let\olddefinition\definition
\renewcommand{\definition}{\olddefinition\normalfont}

\title{The flecnode polynomial: a central object in incidence geometry}
\author{ Nets Hawk Katz  
   \thanks{The author was partially supported by NSF grant DMS 1266104} }

\begin{document}

\maketitle

\section{Introduction}

Incidence geometry as we shall understand it in this lecture is the study of configuration of points and lines in real Euclidean
space, largely with a view towards bounding the number of incidences, that is pairs of points and lines where the point lies on the
line. The study of incidence geometry has a long history. One of its high points is the Szemer\'edi-Trotter theorem
\cite{ST}:

\bigskip

\begin{theorem} \label{SzT} A set of $m$ distinct lines and $n$ distinct points in the Euclidean plane has at most
$$O(n^{{2 \over 3}} m^{{2 \over 3}} + n + m)$$
incidences.
\end{theorem}

\bigskip

One thing that is remarkable about this result, published in 1983, is that except for constants, it is entirely sharp. More delicate
incidence questions, for instance those involving incidences between configurations of points and lines which were in some way forced
to be higher dimensional, for a long time eluded sharp treatments. The subject has undergone a revolution recently, however, in which
sharp results for many problems became possible. The revolution was started by the following result of Ze'ev Dvir in 2008. 
\cite{D}
(Strictly speaking,
in this lecture we view Dvir's result as outside of incidence theory because it is over finite fields.)

\bigskip

\begin{theorem} \label{Dv} Let $q$ be a power of a  prime and $F_q$ be the finite field of $q$ elements. Let $E \subset F_q^n$ be
a set of points containing a line in every direction. Then
$$|E| \gtrsim_n  q^n.$$
\end{theorem}

\bigskip

Here, absolute values denote the cardinality of sets of points, and the notation $\gtrsim_n$ means that we may be concealing a constant
depending on $n$, but certainly not on $q$. Dvir's proof was obtained by studying a polynomial vanishing on the set of points. It easily solved
in the affirmative a conjecture of Tom Wolff's that had been considered quite hard. (It was the finite field Kakeya problem and was intended
as an analog for the real Kakeya problem which arises in geometric measure theory and harmonic analysis.)
Previously the behavior of polynomials had not been used
much in incidence geometry and a number of breakthroughs occurred once it was realized that this is possible.

\bigskip

A lot was accomplished in two papers by myself and Larry Guth. \cite{GK1}, \cite{GK2}
We settled the Joints conjecture in ${\bf R}^3$. Three lines are said to
form a joint at a point $p$ if the lines are all incident to $p$ and not coplanar.

\bigskip

\begin{theorem} \label{GK1}(Joints problem (Guth-K.) A set of $N$ lines in ${\bf R}^3$ forms at most $O(N^{{3 \over 2}})$ distinct joints.
\end{theorem}

\bigskip

We settled a conjecture of Bourgain intended to serve as an analog in incidence geometry for the Kakeya problem.

\bigskip

\begin{theorem} \label{GK2} (Bourgain problem (Guth-K.))Let $E$ be a set of points in ${\bf R}^3$ and let $L$ be a set of $N^2$ lines so that no more than
$N$ lines of $L$ lie in any common plane and so that each line of $L$ is incident to at least $N$ points of $E$. Then
$$|E| \gtrsim N^3.$$
\end{theorem}

\bigskip

We obtained a near-solution to a famous problem of Erd\"os on distinct distances between points in the plane.

\bigskip

\begin{theorem} \label{GK3} (Erd{\H o}s distance result (Guth-K.)) Let $E$ be a set of $N$ points in ${\bf R}^2$. Let $D$ be the set of distances between any two
of the points of $E$ then
$$|D| \gtrsim {N \over \log N}.$$
\end{theorem}

\bigskip

We make some remarks on what Theorems \ref{GK1}, \ref{GK2}, and \ref{GK3} have in common. Each of theorems \ref{GK1} and \ref{GK2}
is clearly a result about Euclidean lines in three dimensions. A completely general set of lines in three dimensions might lie in a plane, in which case,
no incidence result better than Theorem \ref{SzT} is possible. That is why each of Theorems \ref{GK1} and \ref{GK2} contain hypotheses preventing
too many lines from lying in a plane. In the case of Theorem \ref{GK1}, this hypothesis is that each triple of lines forming a joint is noncoplanar. In the
case of Theorem \ref{GK2}, things are more explicit. No more than $N$ lines are allowed to lie in a plane. Theorem \ref{GK3} appears to be different.
It doesn't mention lines in three dimensions at all. But the proof is obtained through what's now called the Elekes-Sharir framework. (Blame me. I
named it thus because I learned about it from a particular paper of Elekes and Sharir. \cite{ES}) This framework is a kind of realization of the Erlangen program.
Instead of studying  Theorem \ref{GK3} in the Euclidean plane where it is stated, we study in the groups of rigid motions which is three dimensional. In
fact, it contains Zariski dense open sets which can be viewed as ${\bf R}^3$. It turns out that proving Theorem \ref{GK3} precisely amounts to solving
an incidence problem between points and lines in ${\bf R}^3$ in which the lines are restricted from being too much in a two-dimensional set.

\bigskip

Theorems \ref{GK1} and \ref{GK2} were discovered before Theorem \ref{GK3} and may each be viewed as special cases of the incidence
result underlying Theorem \ref{GK3}. Their proofs provided essential clues for discovering the proof of Theorem \ref{GK3}. It is hard
to imagine the bigger result coming first. In the last few years, incidence geometry has become much more crowded and many people
are working on and refining the new polynomial methods. The proof of Theorem \ref{GK1} has been so simplified that one can obtain it
without even using Bezout's lemma. But in some sense, these simplifications are merely cosmetic and serve to obscure the unity of the three
theorems. The proof of Theorem \ref{GK3} has essentially two parts. One part is topological. Roughly it serves to show that if the set of points
we are investigating does not lie in the zero set of too low degree a polynomial, one gets a kind of three dimensional improvement of Theorem
\ref{SzT}. While our paper was the first to really do this so that three dimensionality is expressed in terms of polynomials, and this method
is now referred to as polynomial partitioning, it is very much in the spirit of the pre-existing theory of incidence geometry and of the
decomposition method which provided the original proof of Theorem \ref{SzT}.  The second part of the proof is about what to do when
the points in question are in the zero set of a fairly low degree polynomial. Both Theorems \ref{GK1} and \ref{GK2} are special cases of
this part. They are sufficiently easy special cases that they can be resolved quite simply. But what all results in point line incidence theory in
${\bf R}^3$ which we can do because the points lie in the zero set of a low degree polynomial have in common is this: if there are too many incidences,
then most of the offending lines lie in an algebraic ruled surface (of fairly low degree). It is this commonality which I hope to emphasize in the current
lecture.

\bigskip

A number of criticisms can be levelled against the current lecture. The most serious is that it doesn't contain any actual proofs. However,
it does contain sketches which can be turned into actual proofs, at least by me, and which express the way I thought about the results
with Guth when I was working with him on them. A number of standard uniformity arguments are entirely sloughed over and the algebraic arguments
contain slightly excessive assumptions of genericity which have to be justified. Nonetheless, for that perfect reader who catches the zeitgeist, this is supposed to provide a short sweet introduction to the subject, emphasizing major ideas and removing annoying details. A few open problems
are mentioned where they're related to the subject of the lecture. I can't claim any originality in posing them. A lot of ideas about the frontier
were in the air at an IPAM program in Spring 2014. I thank profusely any participant I may have inadvertently stolen from.

\section{The Cayley Salmon theorem }    

In this section, we prove the main result which allows the theory of ruled surfaces to enter incidence geometry. This is
the theorem of Cayley and Salmon which says that any algebraic surface in ${\bf C}^3$ contain enough lines must have
a ruled component. More precisely it says:

\bigskip

\begin{theorem} \label{CS}(Cayley-Salmon theorem) Let $p(x,y,z)$ be a polynomial of degree $d$ on ${\bf C}^3$. Then there is a polynomial
$Flec(p)(x,y,z)$ of degree no more than $11d-24$ which vanishes at a point $w=(x,y,z)$ of the zero set of $p$ only if there is a line
containing $w$ so that $p$ restricted to the line vanishes to third order at $w$. If $Flec(p)$ vanishes at all points of the zero set of $p$
then the zero set is ruled. (That is, through each point of the zero set, there is a line contained in the zero set.) \end{theorem}

\bigskip

Recently, there has been a lot of confusion about the Cayley Salmon theorem. As often happens when people are confused, Terry Tao in the goodness of his heart, posted
an elementary proof of the theorem on his blog, to much acclaim \cite{T}. One thing a bit odd about this is that Salmon also posted an elementary proof. This happened in 1862  \cite{S}, so
he put his proof in an analytic geometry textbook. This book is now past the expiration of its copyright, but unfortunately amazon still charges around twenty bucks, to
produce a copy by print-on-demand. This seems off-putting to people. The Tao blog admits that the theorem dates to at least 1915. 
(Probably this date comes from the most common reprint of the 5th edition.)
One of the motives for this lecture is to defend Salmon's honor and explain his original proof.

\bigskip

There are a number of reasons why Salmon's proof is difficult to parse for modern readers. One is that Salmon's notation is not so good as Tao's. Another is that
Salmon was not restricting his attention to ruled surfaces. He was interested in having similar statements for surfaces ruled by other
classes of curves than arbitrary lines. The idea was that he
assumed his surface had the desired property. Then he derived a differential equation which the surface had to satisfy. Then he
observed that this differential equation had first integrals and that these imply that any surface satisfying his equation also has the desired ruling.
\bigskip

We begin with the first step. We follow the notation of  article 437 of Salmon's book where the argument is explained. We assume that a surface is ruled. Then locally (at least away from singular points), it
can be written as a one parameter family of lines. We write the equations of these lines
$$z=c_1 x+ c_3, \quad y=c_2 x + c_4.$$
We view $c_2,c_3,$ and $c_4$ as arbitrary functions of $c_1$. Of course, this doesn't work entirely in general. The projection of our family of lines into the $xz$ plane might have constant slope. But we change
coordinates so that this is true locally. Then we view the surface as being locally a graph, 
$$z=f(x,y).$$
Of course, a change of variables may be required and $f$ is not a polynomial. It is an algebraic function obtained implicitly from the equation
$$p(x,y,z)=0.$$
At a regular point, we can use the implicit function theorem to solve this for some choice of the $z$-direction. Now, however, we find a third order partial differential equation satisfied by $f$ from the
parametrized description of our surface as a family of lines.

\bigskip

We shall be concerned for the moment with the behavior of $f$ along a single line.
Like Salmon, we adopt the traditional notation for partial derivatives. We let the first derivatives be
$$p={\partial f \over \partial x}(x,c_2 x+ c_4); \quad  q={\partial f \over \partial y}(x,c_2 x + c_4).$$
Here we've emphasized that these partials are being evaluated along one line, where the line is parametrized by $x$. We omit this same dependence on the variables in describing the rest of Salmon's notation.
We let the second derivatives be
$$r={\partial^2 f \over \partial x ^2}; \quad s={\partial^2 f \over \partial x \partial y}; \quad t={\partial^2 f \over \partial y^2},$$
and finally
$$\alpha={\partial^3 f \over \partial x^3}; \quad \beta={\partial^3 f \over \partial x^2 \partial y}; \quad \gamma={\partial^3 f \over \partial x \partial y^2}; \quad \delta={\partial^3 f \over \partial y^3}.$$

\bigskip

Now we write down what it means that an individual line lies in our surface. We have
$$c_1 x + c_3=f(x,c_2 x + c_4).$$
We differentiate this equation in $x$, in effect differentiating along the line and we obtain Salmon's equations
$$p+mq=c_1; \quad m=c_2,$$
from the chain rule. Note that the second equation is in place to say that $m$ is constant along the line. Thus we are free to keep differentiating along the line
as long as we make sure that we follow the chain rule, producing a $c_2$, now called $m$, every time we introduce a partial of $f$ with respect to $y$.
A second derivative produces
$$r+2sm+ tm^2=0.$$
This is not yet a differential equation for $f$ because it still involves $m$, one of the constants of the line. But we take a third derivative:
$$\alpha + 3\beta m + 3 \gamma m^2 + \delta m^3=0.$$
We solve the quadratic equation, obtain two solutions
$$m={-2s \pm \sqrt{4s^2 -4rt} \over 2r}.$$
Plugging each value of $m$ that we obtain  into the cubic equation and multiplying the two equations together, we get an equation which is rational in 
$\alpha,\beta,\gamma,\delta,r,s,$ and $t$. This is the PDE which we assert guarantees that a surface is ruled. (We don't have to know the constants $c_1,c_2,c_3,c_4$
in order to interpret the differential equation. All the partials are evaluated at the same place.) What the equation says precisely is that one of the two complex directions in which the
quadratic form associated to the second derivative vanishes also annihilates the third derivative. In other words, over the complex numbers, the surface has a line tangent to it at third
order. Now what remains is for us to see is that the fact that this equation is satisfied actually implies that the surface contains a line at its generic point.

\bigskip

So once again, we have a surface $z=f(x,y)$ which satisfies our differential equation. Rather than write it out in all its horror, we go backwards by a reversible step and observe that at a generic point, we
have a once differentiable function $m(x,y)$ satisfying
\begin{equation}   \label{second} r+2sm+ tm^2=0, \end{equation}
and
\begin{equation} \label{third} \alpha + 3\beta m + 3 \gamma m^2 + \delta m^3=0, \end{equation}
where now the derivatives $r,s,t,\alpha.\beta,\gamma,\delta$ are viewed as being evaluated on $x,y$ rather than on a line. In the case of planes and quadrics,
the function $m$ can be found by hand. In all other cases, it is produced for us because we have a unique line vanishing to third order at the generic point. Now as before,
we can parametrize the line at a given  $(x,y)$ by $c_1,c_2,c_3,c_4,$ where always $c_2=m$. It is enough to show that $c_1,c_3,$ and $c_4$ are also just functions of $m$. If this is the case, then
it is the same line tangent line vanishing to third order on the points of each level set of $m$ on the surface which implies that the level sets are in fact contained in the lines and that the surface is ruled
by these lines. Now, we can easily write down each of $c_1$, $c_3$, and $c_4$ in terms of $x,y$ and $m$. Namely
$$c_1=p+mq,$$
$$c_3=z-c_1 x,$$
and
$$c_4=y-mx.$$
Our goal now is simply to show each of these three functions is a function of $m$. We will do this by showing that the gradient of each one is a multiple of the gradient of $m$. Thus the level curves are also
level curves of $c_1$, $c_3$, and $c_4$.

\bigskip

We begin with a preliminary calculation. We will differentiate equation (\ref{second}) first with respect to $x$ and then with respect to $y$. With respect to $x$, we get
$$\alpha + 2\beta m + \gamma m^2 + 2s {\partial m \over \partial x} + 2 t m {\partial m \over \partial x} =0.$$
Then differentiating with respect to $y$, we get
$$\beta + 2\gamma m + \delta m^2 + 2s {\partial m \over \partial y} + 2 t m {\partial m \over \partial y} = 0.$$
Adding the first equation to $m$ multiplied by the second equation and dividing by $2s+2tm$
[we leave as an exercise to the reader to work the excluded case where $2s+2tm$ is identically zero], we obtain
$${\partial m \over \partial x} + m {\partial m \over \partial y} = 0.$$

\bigskip

Now we compare the gradients of $c_1$, $c_3$, and $c_4$.
We calculate
$$\nabla c_1=(\alpha + m \beta + {\partial m \over \partial x} q, \beta + m \gamma + {\partial m \over \partial y}  q).$$
Taking the dot product of this with $(1,m)$ and using equation (\ref{second}), we say that $\nabla c_1$ points in the same
direction as $\nabla c_2$.
Further we calculate
$$\nabla c_3=(p-c_1,q) - x \nabla c_1.$$
By dotting with $(1,m)$, we see, using the fact that $\nabla c_1$ is already in the direction of $\nabla m$ and using the definition of $c_1$ as $p+mq$,
we see also that $\nabla c_3$ is in the same direction as $\nabla m$. Finally, we calculate
$$\nabla c_4= (-m,1) - x \nabla m,$$
which is immediately seen to be in the direction of $\nabla m$. This, in effect, is Salmon's argument. He refers to the equations
$$c_1=\psi(m); \quad c_3=\phi(m);  \quad c_4=\chi(m),$$
with $\psi,\phi, $ and $\chi$ as unknown functions of $m$ as the first integrals of his differential equations for surfaces. Part of the reason this proof of Salmon's
is difficult to parse is that he claims it in much greater generality for any surface  ruled by curves of constant complexity. Basically, if
the curves come from a family with a fixed number of parameters and we assume that all parameters are a function of one of the parameters
as we did for lines
 then the differential equation obtained by reducing away all parameters of the curves must imply that the surface is ruled by such curves.

\bigskip

Once this is done, arriving at Theorem (\ref{CS}) is merely a matter of keeping track in the case of a surface $p(x,y,z)=0$ of the
polynomial obtained from checking whether the vectors in the tangent space to a point in whose direction lines vanish to second order
have the property that these lines actually vanish to third order. From the point of view of reduction theory, this is precisely analogous
to the process of eliminating $m$ from the equations (\ref{second}) and (\ref{third}) which we have discussed. This yields a polynomial
of degree $11d-24$ as discussed in Article 588 of Salmon.

\section{On intersections between lines}

\bigskip

In this section, we describe the applications of Theorem \ref{CS} to real incidence geometry. We remark that it is easy to express surprise
that the theorem is applicable at all. After all, the theorem is stated over the complex numbers. Still the reals are a subfield and it is
possible to exploit this. An important and basic tool is the following variant of Bezout's lemma.

\bigskip

\begin{lemma} \label{Bezout} Let $p(x,y,z)$ and $q(x,y,z)$ be two complex polynomials of three variables of degree $m$ and $n$
respectively. Suppose that $p$ and $q$ vanish simultaneously on more than $mn$ complex lines. Then $p$ and $q$ have a nontrivial common factor.
If $p$ and $q$ are both real, then their common factor must be real. \end{lemma}

\bigskip

It may be viewed as a drawback of the flecnode polynomial for investigating real geometry that it is possible to find real polynomials $p(x,y,z)$
for which $Flec(p)(x,y,z)=0$ but nonetheless the real surface $p(x,y,z)=0$ contains no lines. An obvious example is the unit two sphere given
by 
$$p(x,y,z)=x^2+y^2+z^2-1.$$
When we view the zero set of $p$ as a complex surface, it is ruled (and in fact doubly ruled), but over the reals it contains no lines. However, this
is not the way that we ever use Theorem \ref{CS}.  We don't assert that a surface contains many lines by showing its flecnode polynomial vanishes
identically. Instead, we start with a surface containing many lines and conclude that it has a ruling. Indeed when the lines are real, it is often possible to
show that the ruling is real. But it isn't really important. We are interested in intersections between these lines that we already know about and the presence
of the ruling allows us to show that some lines don't intersect, even if the ruling is complex.

\bigskip

To wit, we state the following corollary of Theorem \ref{CS} and Lemma \ref{Bezout}

\bigskip

\begin{corollary} \label{ruled} Let $p(x,y,z)$ be an irreducible polynomial of degree $d$. Suppose the surface
$p(x,y,z)=0$ contains more than $11d^2-24d$ complex lines. Then the surface must be ruled over the complex numbers. \end{corollary}

\bigskip

The proof of Corollary \ref{ruled} is simple. If a line $l$ is in the zero set of $p$, then at each point of the line $l$, there is
a line going through the point namely $l$ on which $p$ vanishes to order at least three. Thus $l$ is in the zero set of $Flec(p)$.
Applying Lemma \ref{Bezout}, we conclude that $p$ and $Flec(p)$ have a nontrivial common factor, and since $p$ is assumed irreducible,
it must be that $p$ is that factor. Thus $Flec(p)$ vanishes on the zero set of $p$ and we conclude from Theorem \ref{CS} that
the zero set of $p$ is ruled over the complex numbers.

\bigskip

This raises something of an open problem. (It is probably not a very serious one.) The corollary above is written in a form that is rather
easily usable by incidence geometers. If an irreducible algebraic surface of low degree contains too many lines then it is ruled. It might be
useful to have such results for other curves and in higher dimensions. The result for curves inside surfaces in ${\bf R}^3$ is probably
already contained in Salmon's Article 431. In general, if one finds a polynomial of sufficiently low degree in ${\bf R}^n$ whose zero set
contains enough $l$ dimensional surfaces of a certain class, does this imply that many of those $l$ dimensional surfaces lie
in a surface of dimension $l+1$ or greater ruled by the $l$ dimensional surfaces. [We can't require something of dimension
greater than $l+1$ because an $l+1$ dimensional surface ruled by the $l$ dimensional ones already contains infinitely many.]
A number of special cases are in the literature (see e.g.  \cite{G2}, \cite{Land})
but maybe somebody who is good at calculus should write a general theorem and
greatly demystify the subject. A fun exercise might be to see whether the higher dimensional joints problem is related to ruled surfaces
in the way we're about to show the regular joints problem is.

\bigskip

In order to utilize Corollary \ref{ruled}, we should ask how can we find a low degree polynomial that vanishes on a set of lines. One
approach is simply to use surface-fitting.

\bigskip

\begin{lemma} \label{fitting} (curve-fitting)  Let $Q$ be a set of $N^3$ points in ${\bf R}^3$. Then there is a polynomial
of degree $O(N)$ vanishing on the points of $Q$. Let $L$ be a set of $N^2$ lines then there is a polynomial of
degree $O(N)$ vanishing on all the lines of $L$. \end{lemma}

\bigskip

The proof of the first part of the Lemma is just that the general polynomial of degree $O(N)$ has more than $N^3$
coefficients. The system of linear equations on the coefficients which says that the polynomial $p$ vanishes on all the points
of $Q$ is underdetermined. To prove the second part, just pick $KN$ points on each line, where $K$ is constant which is large
compared to the implicit constant in the $O(N)$ of the first part. Now a polynomial which has degree $O(K^{{1 \over 3}} N)$
vanishes on all these points. Since $O(K^{{1 \over 3}} N)$ is smaller than $KN$ by the fundamental theorem of algebra, the 
polynomial must vanish on all the lines.

\bigskip

To apply Corollary \ref{ruled}, it should be clear that Lemma \ref{fitting} is useless. The reason is that this fitting applies to
all sets of lines, whereas we are trying to find structure in a set of lines. Luckily we have a technique for finding lower degree polynomials
that vanish on sets of lines when those lines have unusually many intersections.

\bigskip

\begin{lemma} \label{DR} (Degree reduction) Let $L_1$ and $L_2$ be sets of at most $N$ lines. Suppose each line $l$ of $L_2$ intersects at least $QN^{{1 \over 2}}$
lines $l^{\prime}$ of $L_1$ with $Q>0$ a large real number. Then there is a polynomial of degree $O({N^{{1 \over 2}} \over Q})$
which vanishes on all the lines of $L_2$. \end{lemma}

\bigskip

To prove this, we make a random selection $L_3$ of lines from $L_1$ so that each line is chosen independently with probability $\sim {1 \over Q^2}$.
Then with high probability there are $\sim {N \over Q^2}$ lines of $L_3$. Moreover with high probability, each line of $L_2$ intersects
$\sim {N^{{1 \over 2}} \over Q}$ lines of $L_3$. But there is a polynomial $p$ of degree $\sim {N^{{1 \over 2}} \over Q}$ which vanishes on
all the lines of $L_2$. The reader may check that by setting the constants correctly, we can make the degree of the polynomial slightly lower
than the number of lines of $L_3$ each line of $L_2$ intersects. Thus all the lines of $L_2$ are in the zero set of $p$.

\bigskip

A version of the Lemma above was first used in the proof of Theorem \ref {GK1}. These days, people gleefully tell me that no one ever
uses degree reduction to prove the joints theorem. There are simplifications. (See \cite{EKS} ). Isn't it far better just to say, ``let us consider the polynomial
of lowest degree vanishing on all significant lines" and not to worry at all about what that degree is. But it is a remarkable fact that the joints theorem
only works (to within a constant) in the regime where we have significant degree reduction, that is where most of the lines are arranged in 
ruled surfaces. Similarly, we didn't even need to use degree reduction to prove Theorem \ref{GK2}. If the set $E$ contains only ${N^3 \over Q}$
points, just curve-fitting guarantees that the lines are in the zero set of a polynomial of degree ${N \over Q^{{1 \over 3}}}$ and with $Q$ sufficiently
large, this already guarantees that the lines are mostly arranged in ruled surfaces. Why should we care that the degree is really ${N \over Q}$.
(See \cite{G1} for a partial answer.)

\bigskip

We used Lemma \ref{DR} in conjunction with Corollary \ref{ruled} in proving the following result which played a role in the proof of Theorem
\ref{GK3}.

\bigskip

\begin{theorem} \label{GK4} Let $L$ be a set of $N^2$ lines. Suppose at most $O(N)$ of the lines of $L$ lie in a common plane
and that at most $O(N)$ lines lie in a common doubly ruled surface (parabolic hyperboloid or regulus). Then letting $P$ be the set of points contained in at least two
lines then $|P|=O(N^3)$. \end{theorem}

\bigskip

We briefly sketch the proof of Theorem \ref{GK4}. Assume that $|P|=QN^3$ with $Q$ large. The worst case is that there are $\sim N^2$
lines each meeting $QN$ lines. (Situations where the intersections are concentrated on fewer lines end up being easier to handle since with a work
this corresponds to having fewer lines account for all the intersections.)
Then using degree reduction, we find a polynomial $p$ of degree $O({N \over Q})$ which vanishes on these $N^2$ lines. We factor $p$
into irreducible components (over the complex numbers) $p_1 \dots p_d$. Each line is in the zero set of one of the components. Each zero-set
having its share of lines is, in fact, ruled. If a component is a plane or regulus, it has $O(N)$ lines, less than its share. We end up concluding
that most lines lie in ruled components, and as before, if we show these lines aren't involved in most of the intersections, then we end up
with an easier problem.

\bigskip

A line not lying in some ruled surface of degree $k$, will only intersect that ruled surface in at most $k$ points, so it emerges that most
of the intersections we have to worry about come from within an irreducible ruled surface. How often can lines in a
non-planar, non-regulus ruled surface intersect? We say that a point in an irreducible ruled surface is exceptional if it intersects an infinite
number of lines contained in the surface. We say a line in the ruled surface is exceptional if it meets an infinite number of lines contained in the
surface. An irreducible ruled surface contains at most one exceptional point and two exceptional lines. Thus we have at most $O(N)$ exceptional
lines in the whole story, which contribute at most $O(N^3)$ intersections which is harmless. Nonexceptional lines in a ruled surface
of degree $d$ (in Salmon's language: generators) meet exactly $d-2$ other generators of the surface. Again, this gives
every line at most $O(N)$ intersections. We conclude that there cannot be more than $O(N^3)$ points of intersection.

\bigskip

The same ideas prove Theorem \ref{GK1}. To prove the joints theorem, we must take account of the fact that we have
removed the restriction on the number of lines in a regulus and the number of lines in a plane. However a joint cannot come exclusively
from the lines in a plane or from the lines in a regulus. Each line outside a plane or regulus can intersect the plane or regulus at most once or
twice respectively. The bounds on internal intersections in other ruled components control the number of internal joints there.

\bigskip

At the time we proved Theorems \ref{GK1} and \ref{GK2}, we were not aware of the Cayley-Salmon polynomial $Flec(p)$. We thought
these theorems revolved around the gradient and second fundamental forms respectively. But this was short-sighted on our part. Regardless
of which polynomial we use, these results are dramatic only in the ruled regime. This creates rather big problems in operating just outside these
regimes. A question frequently posed by Guth (see \cite{G2}  ) asks: suppose we have a set of $N$ lines making almost $N^{{3 \over 2}}$ joints.
Must a large set of lines (say of size almost $N^{{1 \over 2}}$ be coplanar? We have no good algebraic way of addressing this question
yet because we are just outside the range where degree reduction works. We know no special algebraic properties of the lines.

\bigskip

\section{Elekes-Sharir framework and polynomial partitioning}

\bigskip

Theorem \ref{GK4} plays a role in the proof of theorem \ref{GK3}, but is not the whole of the proof. We now briefly review
the Elekes-Sharir framework which explains how incidences between points and lines in space give information about distances between
points in the plane.

\bigskip

Given $E$ a set of $N$ points and $D$ the set of distances, we may define $Q \subset E^4$ to be the set of distance quadruplets,
namely $(e_1,e_2,e_3,e_4)$ is a distance quadruplet if the distance between $e_1$ and $e_2$ is the same as the distance
between $e_3$ and $e_4$. A simple application of the Cauchy Schwarz inequality shows that
$$|Q| \geq {|E|^4 \over |D|}.$$
Thus to prove Theorem \ref{GK3}, it suffices to prove
$$|Q| =O(N^3 \log N).$$

\bigskip

Now the ancients had a more descriptive term for distance quadruplets. They referred to them as pairs of congruent line segments.
And one thing they knew is that whenever two line segments are congruent there is a rigid motion between them. The space of rigid
motions is three dimensional. It consists of all rotations around a center and all translations. If we restrict to non-translations, a good
coordinate system is given by the center of the rotation in Cartesian coordinates, together with the cotangent of half the angle of rotation.

\bigskip

Given two points of $E$, say $e_1$ and $e_3$, we let $l_{e_1 e_3}$ be the set of rigid motions taking $e_1$ to $e_3$. This is a one
dimensional set and in the coordinate system described in the paragraph above, it is a line. We see that $(e_1,e_2,e_3,e_4)$ form
a distance quadruple precisely when $l_{e_1 e_3}$ and $l_{e_2 e_4}$ intersect. We are back to incidence theory. It is fortunate
that with $L$ the set of $N^2$ lines in rigid motion space of the form $l_{e_1 e_3}$, we have no more than $O(N)$ in a regulus
and no more than $O(N)$ in a plane. Then Theorem \ref{GK4} tells us, that there are at most $O(N^3)$ points of intersection of two
lines from $L$. Unfortunately, this is not enough.

\bigskip

Consider a point where $k$ lines of $L$ meet. This contributes $k^2$ distance quadruplets. We need to keep the number of distance
quadruplets below $N^3 \log N$. If too many of the points where at least two lines meet have many lines meeting there then we're sunk.
We let $P_k$ be the set of points where between $k$ and $2k$ lines of $L$ meet. As long as we can prove that
$$|P_k|=O({N^3 \over k^2}),$$
by dyadically decomposing, we obtain the desired bound. However this is a tricky business. There is no purely algebraic argument. (The
estimate doesn't hold in finite fields.) To take care of this, we use some topology.

\bigskip

Using the polynomial ham sandwich theorem, we obtain the following polynomial partitioning lemma which has proved quite useful.

\bigskip

\begin{lemma} \label{PP} Let $F$ be a set of $M$ points in ${\bf R}^3$. Then for any $s$, a power of 2, there is a real polynomial $p(x,y,z)$ of degree $O(s^{{1 \over 3}})$ so that the complement of the zero set of $p$ in ${\bf R}^3$ has at most $s$ connected components
with points of $F$ and each connected component contains at most ${M \over s}$ points.
\end{lemma}

\bigskip

It is important to note that Lemma \ref{PP} does not guarantee us that most of the points aren't in the zero set of the polynomial. As it turns
out, we will be very happy if they are.

\bigskip

We proceed now to sketch a proof that indeed we have the estimate
$$|P_k| = O({N^3 \over k^2}).$$
We suppose not. Then there are ${QN^3 \over k^2}$ such points with $Q$ large. Note that $Q$ is certainly never larger than
$k^2$ simply by Theorem \ref{GK4}. Better {\it a priori} estimates are possible.
 We would like to subdivide this set of points
into components of size at most $k$. By Lemma \ref{PP}, there is a polynomial $p$ of degree ${Q^{{1 \over 3}} N \over k}$
which does this. We divide into two cases. In the first case, most of the points of $P_k$ are in the complement of the
zero set of $p$. In the second case, most of the points of $P_k$ are on the zero set.

\bigskip

In the first case, we obtain upper and lower bounds on the number of incidences between components and lines. We say a line and a component
are incident, if there is a point on the interior of the component which lies on the line. Let $I$ be the number of such incidences. Since we are
in the case where most points are in the interior, many components have $\sim k$ points. Each point has $~k$ lines through it. Since
any two points have at most one line in common, there is not too much double counting and we get $\sim k^2$ lines incident to each cell.
We conclude
$$I \gtrsim {QN^3 \over k}.$$

\bigskip

On the other hand, each line can switch components only by crossing the zero set of the polynomial. Each line does this at most as many
times as the degree of $p$ plus one. Since there are only $N^2$ lines, we conclude
$$I \lesssim {Q^{{1 \over 3}} N^3 \over k}.$$
We have arrived at a contradiction.

\bigskip

Thus, we are in the second case. Most of the points are in the zero set of the polynomial. There are ${QN^3 \over k^2}$ of these points
each incident to $k$ lines. This gives, on average, ${QN \over k}$ incidences per line. Since the polynomial has
degree ${Q^{{1 \over 3}} N \over k}$, this means that average lines are in the zero set of the polynomial. As usual, the worst case
is that most of the $N^2$ lines are average. (Because if in fact most of the incidence are created by fewer lines, we can basically redo the argument
with a larger $Q$ and smaller $N$.) If we are in the setting where the lines are all average, once again, we are in the domain where the
lines are structured into ruled surfaces, and we can use this structure much as before. Thus we have completed our sketch of the argument
for Theorem \ref{GK3}.

\bigskip

We make a brief remark about the partitioning part of this argument. Contrary to algebra which works best in the complex numbers,
polynomial partitioning seems to work best in the reals. This allows one to prove incidence theorems in the reals much more easily than in the
complex numbers. As an open problem, we suggest considering a set of points in ${\bf C}^3$ no more than half of which is in the
zero set of any low degree polynomial. (For instance, degree lower than ${Q^{{1 \over 3}} N \over k}$, as above.) Can one make arguments
giving results analogous to the case above with most points in the complement of the zero set?

\end{document}